\documentclass[12pt]{article}
\usepackage{amssymb}
\usepackage{latexsym}
\usepackage{exscale}

\renewcommand{\rq}[1]{(\ref{#1})}
\newcommand{\ov}{\overline}

\newcommand{\dsize}{\displaystyle}

\renewcommand{\div}{\mbox{div}}
\newcommand{\R}{{ \bf R  }}

\newcommand{\ra}{\rightarrow}

\newcommand{\BEL}{\begin{equation}\label}
\newcommand{\EE}{\end{equation}}

\renewcommand{\medskip}{\vskip .5 cm}
\newtheorem{Thm}{Theorem}[section]
\newtheorem{Lemma}[Thm]{Lemma}
\newtheorem{Cor}[Thm]{Corollary}
\newtheorem{Prop}[Thm]{Proposition}
\begin{document}

\centerline{\Large \bf Minimal support results for  Schr\"odinger  equations}

\medskip
\bigskip
\centerline{Laura De Carli, Julian Edward, Steve Hudson and Mark Leckband}

\bigskip

\noindent{\bf to be published in  Forum Mathematicum by De Gruyter}

\bigskip

{\centerline{\bf Abstract}} We consider a number of linear and non-linear
boundary value problems involving generalized
Schr\"odinger equations. The model case is $-\Delta u=Vu$ for $u\in W_0^{1,2}(D)$
with $D$ a bounded domain in ${\bf R^n}$. We use the Sobolev
embedding theorem, and in some cases   the Moser-Trudinger inequality and the Hardy-Sobolev inequality, to derive
necessary conditions for the existence of nontrivial solutions.

These conditions usually involve a lower bound for  a product of
powers of the norm of $V$,  the measure of $D$, and  a sharp Sobolev
constant. In most cases, these inequalities
 are best possible.

\medskip\noindent {\bf Mathematics Subject Classification:} Primary 35B05;
Secondary 26D20, 35P15.

\section* {1. Introduction }
\setcounter{section}{1} \setcounter{Thm}{0} \setcounter{equation}{0}

\bigskip

We show that solutions of  certain   second order elliptic differential equations cannot vanish
on the boundary of arbitrarily small domains. Our first example is  the Schr\"odinger equation
 \BEL{Sch} -\Delta u = V(x) u, \qquad x\in D, \EE
where $\Delta = \sum_{j=1}^n  \partial^2_{x_j}$  is the Laplace
operator, and $D$ is, here and throughout the paper, a bounded domain  of $\R^n$.  Unless stated otherwise,  $u \in
W^{1,2}_0(D)$ will be a solution in the distribution sense    of \rq{Sch},
i.e.,
$$ \int_D \nabla u(x) \cdot \nabla\psi(x) dx= \int_D V(x)
u(x)\psi(x)dx $$ for every $\psi\in C^\infty_0(D)$.

We say that a solution $u$ of \rq{Sch} is  {\it
trivial} if $u(x) = 0$ almost everywhere (a.e.) on $D$. Unless
otherwise stated, we assume  $u$ is nontrivial and complex
valued, and that \rq{Sch} holds in the sense of
distributions.

$V$ is often referred to as   the {\it potential} of the equation. Although several of our proofs work
for complex valued $V$, we will assume throughout this paper that $V$ is real-valued, and we  let
$V_+ = \max\{V,\,0\}$. We denote by $||V||_r =||V||_{L^r(D)}$ the usual Lebesgue norm, and say $u$
 is in the Sobolev space
$W^{1,2}_0(D)$ if $u$ is in the closure of $C_0^{\infty}(D)$ with respect to the norm $ \| \nabla u\|_2$.

Let $K_q(D) $  be the operator norm of the Sobolev embedding $W^{1,2}_0(D)\rightarrow L^{ 2q}(D)$. That is,
 \BEL{sob} K_{q}(D)=\sup_{u\ne 0}  \frac{||u||_{2q}}{ ||\nabla u ||_2}.
 \EE
For $n>2$,  let $2\ov q= \frac{2n}{n-2}$ be the critical index in the Sobolev embedding theorem. That is,
$K_{q}(D)$ is finite for $q\le \ov q$ and infinite otherwise. When $n \le 2$, let $\ov q=\infty$.
Note that when $n= 2$ and $q= \ov q = \infty$, then $K_{q}(D)= \infty$, an exceptional
 case addressed carefully in Section 2. When $n=2$, we will assume $q<\infty$ unless specified
 otherwise.

We often assume   that $V\in L^r$, where $r= q^*$ is the H\"older
conjugate exponent to $q$. So $\frac 1q+\frac 1{r}=1$, with the
convention that $1^* = \infty$.
We will abuse notation slightly, and write $r \ge \frac n2$ to mean $\max
\{1,\frac n2\} \leq  r \le +\infty$. Note that when $r \ge \frac
n2$, then $r^* \leq \frac{n}{n-2}= \ov q$.

Our first result is our simplest, and is central to the rest of the
paper.
\begin{Thm}\label{ThmX} Suppose that the Schr\"odinger equation \rq{Sch}
  has a nontrivial solution $u\in W^{1,2}_0(D)$  for some $V\in L^r(D)$, with $r >\frac n2$ and  $q= r^* $. Then,
  \BEL{main} K_{q}^{2}(D)  ||V_+||_r \ge 1.  \EE
Denote by $u_*$ an extremal \footnote{$u_*$ is an extremal for
\rq{sob} if $ K_{q}(D)=   \frac{||u_*||_{2q}}{ ||\nabla u_* ||_2}.$}
  in the  inequality \rq{sob}. With $u=u_*$, and  \BEL{V}
V(x)= \frac{ ||\nabla u_*||_2^2 }{\| u_*\|_{2q}^{2q}}|
u_*(x)|^{2q-2},\EE
equality is  attained in \rq{Sch}  and in \rq{main}.
\end{Thm}
\medskip
\noindent
When $q<\ov q$, $K_q(D)  $ depends on the volume as well as the
shape of $D$. Let $K_q^*$ be the Sobolev constant associated with
the ball of volume 1. It is well-known that $\dsize
K_q^*=\max_{|D|=1}K_q(D)$. A simple dilation argument proves the
following

\begin{Thm}\label{CorX} Under the assumptions of  Theorem \ref{ThmX},
  \BEL{main*} (K_q^*)^{2} |D| ^{\frac 2 n -\frac{1}{r}} ||V_+||_r \ge 1.  \EE
\end{Thm}
\medskip
\noindent
 So, if  the Schr\"odinger equation \rq{Sch}  has
 nontrivial solutions in $W^{1,2}_0(D)$,  and $||V||_r$ is fixed, then $|D|$
cannot be too small; we say that the
solutions have a {\it minimal support property}.

 The proofs of Theorem \ref{ThmX}, and many others in this paper,
 follow a pattern from Theorem 4.1 in   \cite{DH1}, which we refer to as
  the {\it minimal support sequence}.   For example,
  to prove \rq{main} for   $u\in C_0^2(D)$, we use Sobolev's inequality \rq{sob},
  Green's identity, and H\"older's inequality:
\begin{eqnarray}\nonumber ||u||_{2q}^2
& \leq &  K_q^2(D)   \int_D |\nabla u|^2 dx
 =       - K_q^2(D)  \int_D   \ov u \Delta u\ dx \nonumber \\
& = &      K_q^2(D)  \int_D   |u|^2V\ dx
   \leq   K_q^2(D)  \int_D   |u|   ^2V_+ \ dx \nonumber \\
&\leq & K_q^2(D)  ||u^2||_q ||V_+||_r
=  K_q^2(D)  || u||_{2q}^2 ||V_+||_r
\label{minimal}
\end{eqnarray}
which implies \rq{main}. If $u=u_*$  is an extremal for the Sobolev
inequality \rq{sob}, the first inequality in \rq{minimal} is an
equality. We prove in Lemma \ref{EuLa}   that $u_*$ solves \rq{Sch},
with $V=V_+$ as in \rq{V}. Remarkably, this result also makes
H\"older's inequality into an equality. So, \rq{main} is an equality
too.  This apparent coincidence has an explanation; the solutions of
$-\Delta u=V u$ minimize certain energy functionals. Since $u_*$
solves two similar optimization problems, it has certain unexpected
properties; see \cite{H}. We complete the proof of Theorem
\ref{ThmX} in Section 2.1.
\medskip

Theorem \ref{ThmX} can be interpreted as a necessary condition for
zero to be an eigenvalue for the operator $-\Delta-V$ with quadratic
form domain   $W_0^{1,2}(D)$.

\begin{Cor}\label{ThmY}
  Suppose that the differential equation
  $$-\Delta u -Vu=Eu $$ has  a nontrivial solution in $ W^{1,2}_0(D)$ for some constant $E\leq 0$.  Suppose
  $V$, $r$, $q$ and $K_q(D)$ are as in Theorem \ref{ThmX}.  Then  \BEL{mainY} K_q^2(D) ||V_+||_r \ge 1.  \EE If $n \geq 3$ and $r =\frac n2$, then $K_{\ov q}^2  ||V_+||_{\frac n2} > 1.$
\end{Cor}

\medskip
\noindent
This result for $r>\frac{n}{2}$ follows immediately from
Theorem \ref{ThmX} by noting that $(V+E)_+(x)  \leq V_+(x)$; the
case $n\geq 3$ and $r=\frac{n}{2}$ follows similarly from Theorem
\ref{ThmX3}.

\medskip

The analogue of Corollary \ref{ThmY} for $D={\bf R}^n$ has a long
history in the mathematical physics community, motivated by
questions of the existence of bound states, i.e., $L^2$ eigenvalues,
for the Schr\"odinger operator in ${\bf R}^n$. Specifically,
consider
\begin{equation}
-\Delta u -Vu=Eu, \quad u\in W^{1,2}({\bf R}^n),\label{mp}
\end{equation}
with eigenvalue $E<0$.  For fixed $V$, let $\tilde{N}$ be the number
of negative eigenvalues of $-\Delta  -V$.  The following inequality
is due to Cwickel \cite{Cw}, Lieb \cite{L}, and Rozenblum
\cite{Roz}: \BEL{CLR} C_n\| V_+\|_{\frac n2}^{\frac n2 } \geq \tilde{N}, \quad
n\geq 3, \EE  where   $C_{n}$  depends only on n.  For more on the values
of $ C_{n}$, and also for the cases $n=1$ and $n=2$, the reader is
referred to the survey on bound states by Hundertmark \cite{Hun}.


 A more abstract version of the Cwickel-Lieb-Rozenblum inequality,
 which  applies on a bounded domain $D$, as does our Corollary \ref{ThmY}, is
 derived in Theorem 2.1 in  \cite{FLS}.
 It implies
 \begin{equation}
 C ||V_+||_r \ge \tilde{N},\quad r\geq \frac{n}{2}, \label{Lieb}
\end{equation}
 where the best constant $C$ is unknown, but satisfies $ K_q ^{2}(D) \leq C \leq e^{1-\frac 1r} K_q ^{2}(D)$.
If there is exactly one negative eigenvalue, then $\tilde{N}=1$,
 and thus $C ||V_+||_r\ge 1$; in this special case, the bound \rq{mainY} improves on \rq{Lieb}.
 It is not clear whether \rq{Lieb} can be compared with \rq{main}, since Theorem \ref{ThmX}  involves
  a zero eigenvalue.
\medskip

We now consider the lower bound $r>\frac n2 $ in Theorems \ref{ThmX} and \ref{CorX}.
When  $n\ge 3$ and $r=\frac n2 $, \rq{sob} still holds,
but equality is not attained on any proper subset  $D \subset R^n$.
So, \rq{main} in Theorem \ref{ThmX} (and likewise \rq{minimal} and \rq{main*})
still holds with the same proof, but equality cannot be attained.
However, the estimate is still sharp; see Theorem \ref{ThmX3} in Section 2.2. These results do not extend to  $r<\frac n2$; see  Theorem \ref{NewCn}.

When $n=1$, the critical case is $r=1$. Theorem \ref{ThmX} is still valid, but to attain equality in \rq{main}, we must allow $V$ to be a finite measure  rather than a $L^1$ function; see Section 4.1.

When  $n=2$, Theorem \ref{NewCn} shows that Theorem \ref{ThmX} does
not extend to $r = \frac n2 =1$. The minimal support sequence fails
because $W^{1,2}_0(D)$ does not embed into $L^{\infty}(D)$. However,
we prove an analogue of Theorem \ref{ThmX} when   $V$ is in the
Orlicz space $L\log L(D)$. Our main result in Section
2.4 uses a norm $\| . \|_{N_D} $ for $L\log L(D)$ defined by
\rq{lux}.
\begin{Thm}\label{Orl} Assume that \rq{Sch} has a nontrivial solution with $V$ in $L\log L(D)$ with $D\subset \R^2$.  Then
\BEL{ThmOrl} \dsize \frac{C_2|D|}{4\pi} \| V_+\|_{N_D}\geq 1, \EE
where $C_2$ is the constant of the Moser-Trudinger inequality
\rq{MT}.

Let $u_*$ be an extremal for \rq{MT},   normalized by $\|
\nabla u_*\|_2=1$.   Then for
\BEL{orlV} V = \frac{e^{4\pi|u_*(x)|^2}}{\int_D |u_*(x)|^2e^{4\pi| u_*(x) |^2}
dx},\EE equality holds in \rq{ThmOrl}.
 \end{Thm}

\medskip
\noindent
We have the following analogue of Corollary \ref{ThmY} in this case.

\begin{Cor}\label{OrlY}
 Let  $u\in W_0^{1,2}(D)$ be a nontrivial solution of
\begin{equation}
-\Delta u -Vu=Eu\label{pq}
\end{equation}
 with $E\leq 0$ and  $V\in L\log(L)$. Then \rq{ThmOrl} holds.
\end{Cor}

\medskip
\noindent
To the best of our knowledge, this result is new. For other
results relating the spectrum to $V$ on bounded domains, see
\cite{Hen}.

\medskip

 Theorem  \ref{ThmX} extends, in part,  the main result in \cite{DH2}; two of the authors proved  that
  if $V \in L^\infty(D)$,
and if $u\in C_0(\ov D)$ is a nontrivial solution of \rq{Sch}, then
\BEL{supp1} |D|^{\frac 2n} \left( j^{-1} \omega_n^{-\frac
1n}\right)^2 \cdot||V|| _\infty \ge 1 , \EE where $j$ is the first
positive zero of the Bessel function $J_{\frac n2-1}$.  Equality is
attained when $u=u_* (x)= |x|^{1-\frac n2} J_{\frac n2-1}(|x|)$. The
proof in \cite{DH2}, which is quite different from the arguments
appearing here, compares the level sets of $u(x)$ and $u_*(x)$.

 By  comparing \rq{supp1} and \rq {main*}, we can see at once that the Sobolev constant $K_1^*$
  associated with the ball of volume $1$ is
$ K_1^*=     (j \omega_n^{\frac 1n})^{-1}. $ Since  $K_1^*$  is the
reciprocal of the first eigenvalue of the Dirichlet Laplacian on the
ball of volume $1$, its value is well known, but it is interesting
to observe how this explicit expression follows from our theorems.

\medskip
 Our Theorem \ref{ThmX}, as well as many other theorems in this
paper, can be viewed as a unique continuation result  for solutions
of the equation under consideration. That is,  if $u$ is a solution
of  \rq{Sch} that vanishes on the boundary of $D$, and if $|D|$ is
too small, then $u\equiv 0$ in $D$. In unique continuation problems
the zero set of $u$ is usually assumed to be an open set, or a
point, but in our case, it may be an  $(n-1)$-dimensional boundary. Our assumption that $V\in L^r(D)$, with $r>\frac n2 $, is also critical in these problems.

\medskip
 This paper  is organized as follows: in Section 2  we prove
necessary conditions for the existence of nontrivial solutions for
the Schr\"odinger equation \rq{Sch} with various assumptions on $V$.
In Section 3, we consider similar questions for other well-known
linear and nonlinear second order equations. Minimal support
problems in $\R^1$ are handled separately in Section 4. We have
collected some technical lemmas, perhaps not entirely new, into an
appendix.

\bigskip

\bigskip

\section*{ 2. The Schr\"odinger equation.}
\setcounter{section}{2} \setcounter{Thm}{0} \setcounter{equation}{0}

\bigskip

In this section we prove necessary conditions for the existence of
nontrivial solutions in $W^{1,2}_0(D)$ of the
Schr\"odinger equation $-\Delta u = V(x)u$. We also consider
potentials which do not necessarily belong to $L^r(D)$ with $r>\frac
n2$, but are dominated by a Hardy potential, or belong to $L^{\frac
n2}(D)$ or to an Orlicz space.

\subsection{Completion of the proof of Theorem \ref{ThmX}}

\medskip

We  used the minimal support sequence in the introduction to prove
that $K_{q}^{2}(D) ||V_+||_r \ge 1$ when $u\in C^2_0 (D)$.  When
$u\in W^{1,2}_0(D)$, we cannot apply the standard  Green's identity,
but   we use instead the identity \rq{greenh2} in Lemma
\ref{GreenG}:
$$\int_D |\nabla u|^2 dx = \int_D   |u|^2V\ dx.$$
All other inequalities in the minimal support sequence  hold also
when $u\in W^{1,2}_0(D)$, and so \rq{minimal} is proved.
\medskip

We now prove that equality can occur in \rq{main}. Since $q=r^* <
\ov q$, there exists an extremal $u_* \ge 0$ for the Sobolev
inequality \rq{sob}; this result is probably known, but we prove it
in the appendix as Lemma \ref{Markab} for completeness. Furthermore,
by Lemma \ref{EuLa}, $u_*$  is a  solution in the distribution sense
of \rq{Sch} with $V(x)=
  ||\nabla u_*||_2^2 |
u_*(x)|^{2q-2} /\| u_*\|_{2q}^{2q}$.
 Since $(q-1)r=q$,
$$\left(\int_D |u_*(x)|^{(2q-2)r}dx\right)^{\frac 1r} = \|u_*\|_{2q}^{2q-2}$$
and
$$ K_q^{2}(D) \| V\|_r =  K_q^{2}(D) \frac {||\nabla u_*||_2^2}{\|u_*\|_{2q}^2 } = 1.\ \   \Box  $$

\bigskip

\subsection{A critical case for $n\ge  3$:  $V\in L^{\frac n2}$}
\bigskip

\medskip
 In this section, we assume $n\ge 3$ and  $V\in L^r(D)$, where $r= {\frac n2}$. This
  assumption on $V$  is weaker than the assumption $r>{\frac n2}$ in Theorem \ref{ThmX}, and $r=\frac n2$
   may be regarded as a critical case;
   see also the next two subsections. In Proposition \ref{NewCn}, we show that no minimal
    support result is possible with smaller $r$ by providing explicit counterexamples.
     We also briefly discuss the case $n=2$ there, with more about that in Sections 2.3 and 2.4.
     For $r={\frac n2}$, we have $r^* = \ov q=\frac{n}{ n-2}$. The Sobolev inequality
\BEL{sob2}||u||_{2\ov q}  <  K_{\ov
q}\|\nabla u\|_2. \EE   is strict (since $D\not = R^n$) and
dilation invariant.  So, $K_{\ov q}$ and the corresponding minimal support sequence are independent of $|D|$.  In a celebrated theorem, Talenti (see \cite{T}) proved that
$  K_{\ov q}= (n(n-2)\pi)^{-\frac 12} \left( \frac{ \Gamma(n)}{
\Gamma(\frac n2) }\right)^{\frac 1n}. $ In the next theorem,
instead of a  minimal support result, we prove a "minimal potential
result".
\medskip
\noindent
\begin{Thm}\label{ThmX3}
Suppose $V\in L^{\frac n2}(D)$, with $n\ge 3$. If \rq{Sch}  has
nontrivial solutions in $W^{1,2}_0(D)$, and $\ov q$ and $K_{\ov q}$ are as above,
\BEL{strict} K_{\ov q}^2 ||V_+||_{\frac n2} >  1\EE and the constant
1 is sharp.
\end{Thm}

\medskip
\noindent {\it Proof.}  The minimal support sequence \rq{minimal}
proves \rq{strict}. In this case, we get a strict inequality because
\rq{sob2} is strict. Now we prove that \rq{strict} does not hold if
1 is replaced by any larger constant. Since  \rq{sob2} is invariant
by dilation, the constant $K_{\ov q}$ is independent of  $D$ (and is
also the same for ${\bf R}^n$).  In our proof we will define a
suitable large disk  $D$ and $u\in W^{1,2}_0( D)$ such that $K_{\ov
q}^2||V||_{\frac n2} \approx 1$. Let $v(\rho) =
(1+\rho^2)^{\frac{2-n}2}$, with $\rho=|x|$. Talenti showed
 in \cite{T} that this $v$ gives equality in \rq{sob2} on $\R^n$.
Define $V_v(\rho)= -\Delta v(\rho) / v(\rho)$.

Recalling that the Laplacian of
 a radial function $u$ in ${\bf R}^n$ is
$ \Delta u = u_{\rho \rho}+\frac{(n-1)}{\rho}u_{\rho}, $ it is easy
to verify that $V_v(\rho)= \frac{(n-2) n}{\left(\rho^2+1\right)^2}$
 and that
\begin{equation}
K_{\ov q}^{n} \int_{{\R}^n}|V_v|^{\frac n2} dx = 1.\label{extremal}
\end{equation}
Indeed,
\begin{eqnarray*}
||V_v||_{\frac n2}^{\frac n2}  & = & (n(n-2))^{\frac n2}|S^{n-1}| \int_0^\infty \rho^{n-1} (1+\rho^2)^{-n}d\rho
\\
&=&
\frac{1}{2}  (n(n-2))^{\frac n2}
|S^{n-1}|\int_0^1 \left(t-t^2\right)^{\frac{n}{2}-1}dt =
 (n(n-2))^{\frac n2} \pi^{\frac n2}
\frac{\Gamma \left(\frac{n}{2}\right)}{\Gamma (n)}= K_{\ov
q}^{-n}.\end{eqnarray*}
We have  used the substitution   $
(1+\rho^2)^{-1} = t$. Let $D' =B_R(0)$, where $R >0 $ will  be
specified later. We  will define a compactly supported,  non-negative
function $u$ by perturbing $v$  on   $D \setminus D'$,
 while keeping $K_{\ov q}^2 \| V_u\|_{\frac n2 }\approx 1$.
   Let  $$ u(\rho)=\cases{
   v(\rho)& if $0\leq \rho\leq R$,\cr
   a\rho+ b & if $R < \rho < R+1$, \cr
   c\rho^{2-n} +d & if    $R+1 \le \rho \le \hat{R}$\cr}
   $$
where $a$, $b$, $c$,  $d$ and $\hat{R}$ are chosen below so that $u$ is differentiable,
 and vanishes at $\rho = \hat R = \partial D$. Note that $u$ is harmonic for $\rho>R+1$.
  Choose $a= v'(R)= (2-n) R \left(R^2+1\right)^{-\frac n2}$,
and $b= v(R)- aR =\left(R^2+1\right)^{-\frac n2} \left(R^2(n-1)+1\right)$, which makes $u$ differentiable at $|x|=R$.

In what follows, $C$ will denote a positive constant that may change
from line to line, but is always independent of $R$. When
$R<\rho<R+1$, $u(\rho)= a\rho+b \geq C R^{2-n}$, and $|\triangle u
(\rho)|= |a|(n-1) \rho^{-1} \leq  C(n-1)(n-2)(1+R^2)^{-\frac n2} \leq C
R^{-n}$. Thus, \BEL{I2}\int_{R\le \rho<R+1} |V_u (x)|^{\frac n2}\ dx
\le C \int_R^{R+1} |R^{-2}|^{\frac n2} \rho^{n-1}d\rho \leq C
R^{-1}. \EE

Next, choose $c=R (R+1)^{n-1} \left(R^2+1\right)^{-\frac n2},$ and
$d=( (1-n)R+1) \left(R^2+1\right)^{-\frac n2}$ so that $u$ is
differentiable at $|x|=\rho =R+1$. Since $d<0$,  there exists  $\hat
R > R+1$ for which $u(\hat R)=0$. We let $D = B_{\hat R}(0)$. Then
$u\in W^{1,2}_0(D)$. Since $u$ is harmonic for $\rho>R+1$,
\rq{extremal} and \rq{I2} imply
$$K_{\ov q}^{n}  \int_D |V_u|^{\frac n2} \ dx \le  1 + \frac C R \to 1 $$
as $R\to \infty$. $\Box$

\medskip

   The following constructions show that the conclusions of Theorems \ref{ThmX}
   and \ref{CorX}  do not hold when $r<\frac n2$, nor when $r=1$ and $n=2$.

\begin{Prop}\label{NewCn} Let $n\ge 3$ and $r<\frac n2$ (or $n=2$ and $r=1$). For every $\epsilon >0$, we can find a
non-negative $V_\epsilon\in L^r(B_1(0))$, and a nontrivial
solution $u \in W^{1,2}_0(B_1(0))$ of  $-\Delta u=V_{\epsilon} u $,
such that $\dsize \lim_{\epsilon\to 0}\|V_{\epsilon}\|_r=0$.
\end{Prop}

\medskip
\noindent {\it Proof.} Suppose $n=3$.    Let $\epsilon >0$ be small.
For $\epsilon\le \rho\le 1$, set $u(\rho)=  \rho^{-1} -1$, so
$u(1)=0$ and $u$ is harmonic. For $0\le \rho \le \epsilon$, set
$u(\rho )= a-b\rho^2$. We choose  $b = (2\epsilon^{3})^{-1}$ so that
$u'(\rho)$ is continuous at $\epsilon$, and we chose $a=
\frac{3}{2\epsilon} -1$, so that $u(\rho)$ is continuous at
$\epsilon$. So, $V_{\epsilon}$ is supported on $B_\epsilon(0)$, and
there $\triangle u = -6b$ and $u \geq C\epsilon^{-1} $ (since
$b\rho^2\le b\epsilon^2= \frac{1}{2\epsilon}$). Hence,
$$
|V_{\epsilon}|  \leq C\epsilon^{-2},$$ so for $r<\frac 32$,
$$||V_{\epsilon}||_r^r \leq C \epsilon^{3-2r} \to 0,\mbox{ as}\ \
\epsilon \to 0 .$$
For larger $n$, we set $u(\rho)=  \rho^{2-n} -1$ for $\rho
>  \epsilon$ instead, with a similar proof.
\medskip

For $n=2$, we set $u(\rho)=  -\ln(\rho)$  for $\epsilon\le \rho\le
1$ which is harmonic. For $\rho <\epsilon$, set $u(\rho )=
a-b\rho^2$. We choose $b = (2\epsilon^2)^{-1} $ so that $u'(\rho)$
is continuous at $\epsilon$, and we chose $a= \frac 12 -
\ln(\epsilon)$, so that $u(\rho)$ is continuous at $\epsilon$. Near
0, $\triangle u = -4b$  and $u> -\ln(\epsilon)$; thus, $  0\le
V_{\epsilon}(r)< \frac{2}{\epsilon^2 \ln( \epsilon^{-1})} , $ and
  $ ||V_\epsilon||_1 \le \frac{C}{\ln( \epsilon^{-1})} \to 0$. $\Box$

\medskip
\noindent {\it Remark:}  When  $n=2$,  the proof of Proposition
\ref{NewCn} shows that $|x|^2 V_\epsilon (x)\le \frac{2}{\ln(\epsilon^{-1})} \to 0$, as
$\epsilon \to 0$. See also the remark following Theorem
\ref{NewHS1}.

\bigskip

\subsection{ Hardy potentials}
\bigskip

 We now prove  minimal potential results for solutions  of the  Schr\"odinger equation
 with pointwise bounds on $|V|$, but no longer assuming $V\in L^{\frac n2}(D)$.
 For example, we study $V=C|x|^{-2}$, which is known as a {\it Hardy potential}.
 Let ${\rm  dist}(x) =\inf \{|x-y|, y\in\partial D\}$.

\medskip
\begin{Thm}\label{NewHS1} Suppose $n\ge 2$, and that a measurable $V$
satisfies one of these on $D$:
\begin{itemize}
\item[i)]  \BEL{Hardy}   |V(x)| \leq
\left(\frac{n-2}{2}\right)^2  |x|^{-2}; \quad  or\EE

\item [ii)] $D$ is convex with  piecewise-smooth boundary, and
 \BEL{dist}  |V(x)| \leq \frac 14\ {\rm  dist}(x)^{-2}.\EE
 \end{itemize}
 Then \rq{Sch} has no nontrivial solutions in  $W^{1,2}_0(D)$.
\end{Thm}
\medskip
We do not assume that $0\in D$. Also, note that when $n=2$,
\rq{Hardy} reiterates that \rq{Sch}  has only trivial solutions when $V\equiv 0$.

 When $n=2$ and \rq{Hardy}  is replaced by   $|V(x)| \leq C|x|^{-2}$, for some positive constant $C$, the remark after the proof of Proposition \ref{NewCn} shows that \rq{Sch} can have nontrivial solutions.
\medskip
\noindent {\it Proof.}  Suppose $V$ satisfies \rq{Hardy} and that \rq{Sch} has a
nontrivial solution $u\in W^{1,2}_0(D)$. We use a variation of the classical Hardy Sobolev inequality (see \cite{BV}):
\BEL{H1}\int_D|\nabla u(x)|^2 dx-\left(\frac{n-2}{2}\right)^2\int_D
\frac{|u (x)|^2}{|x|^2}dx\ge C(D)
 \| u\|_2^2 > 0.
 \EE
By Green's identity (Lemma \ref{GreenG}, with $a\equiv b\equiv 1$)
and the above,
\begin{eqnarray*}
\left(\frac{n-2}{2}\right)^2\int_D \frac{|u (x)|^2}{|x|^2}\ dx & < &
\int_D|\nabla u(x)|^2 dx \\ &=&   \int_D V(x)|u(x)|^2 \, dx \leq
\int_D |V(x)||u(x)|^2\,dx
\end{eqnarray*}
which contradicts \rq{Hardy}.

\medskip
Assume now that $V$ satisfies \rq{dist} and that   $u$ is a
nontrivial solution of  \rq{Sch}.  Using an inequality  in \cite{BFT},
  \BEL{BM}
\int_D|\nabla u(x)|^2 dx-\frac 14 \int_D \frac{u^2(x)}{{\rm
dist}(x)^2}dx\ge \frac{1}{4\,\mbox{diam}^2(D)} \int_D |u(x)|^2 dx
>0. \EE By Lemma \ref{GreenG} and the above,
\begin{eqnarray*}  \frac 14 \int_D \frac{|u(x)|^2}{{\rm dist}(x)^2}dx  &<& \int_D|\nabla u(x)|^2 dx \\ &=&
 \int_D V|u(x)|^2\,dx \leq \int_D |V (x)|\,|u(x)|^2\,dx
 \end{eqnarray*}
contradicting \rq{dist}. $\Box$

\medskip
\noindent{\it Remark.}   In the recent paper \cite{FL}, the
authors improve the inequality \rq{BM}, and prove a
Cwikel-Lieb-Rozenblum type inequality for the
 negative eigenvalues of $H=-\Delta - (2D_\Omega)^{-2}+V$,
where $D_\Omega$ is a function that can be replaced by ${\rm
dist}(x)$ when $D$ is convex. They also observe, using a  minimal
support sequence similar  to ours, that $H$ will have no
negative eigenvalues  if $||V_-||_{\frac n2}$  is sufficiently small. For other results
 related to  the first and second parts of our theorem, see
 \cite{Da},\cite{KO}.

\bigskip

\subsection{A critical case for n=2: $V \in L$ log $L$.}

 In this section, we   prove Theorem \ref{Orl} by proving the equivalent
Theorem \ref{Orl2} below. We have observed that Theorem \ref{ThmX}
does not hold when $n=2$ and $V\in L^1(D)$. Here, we  assume $V$ in
the Orlicz space $L$log$L(D)$, defined as  the set of measurable
functions $f$ such that $\int_D |f|(1+ \log^{+} |f|)dx $ is finite.
We will use the Moser-Trudinger inequality (see \cite{M}) as a
substitute for \rq{sob}; it is
       \BEL{MT} \int_D \left(  e^{ 4\pi \left (\frac{|u(x)|}{\|\nabla u  \|_2}\right)^2}
        - 1\right )\  dx \leq C_{2}|D|,\quad  u \in W_0^{1,2}(D), \EE
        where the constant $C_2$ does not depend on $u$ or $D$.  Let $M(x) = e^x -1$, and
 \BEL{N}
 N(y) = \cases{ y \log(y)-y+1, & if $y\geq 1$,\cr  0 &otherwise.\ \cr}\EE
Following \cite{KR}, we set for $V\in L\log L(D)$
\BEL{lux}\|V\|_{N_D} = \inf \left\{ \lambda +\frac{\lambda}{C_2|D|}
\int_D N\left(\frac{|V(x)|}{\lambda}\right) dx;\ \lambda
>0 \right\} <\infty.\EE One can verify that $\|\cdot\|_{N_D}$ defines a norm.
For fixed $V$, we set $F(\lambda) = \lambda
\int_D N\left(\frac{V_+(x)}{\lambda}\right) dx$, so that $\|V_+\|_{N_D} = \inf \left\{ \lambda + \frac{F(\lambda)}{C_2|D|}\right\}$.
\begin{Thm}\label{Orl2}
Suppose that  \rq{Sch} has a nontrivial solution for $V\in L \log
L(D)$, where $D\subset \R^2$. Then, for every $\lambda>0$,
\BEL{inf}  \lambda C_{2} |D| + F(\lambda) \geq 4\pi.
  \EE
Equality can be attained  in \rq{inf} when $u_*$ is an extremal for
\rq{MT}
 and $V=V_+$ is as in \rq{orlV}.
\end{Thm}
 \medskip\noindent\
Theorem \ref{Orl} follows immediately.

\medskip
\noindent
{\it Proof of Theorem \ref{Orl2}.} Fix $u,V$. Let $U =4\pi
\frac{|u(x)|^2}{\| \nabla u\|_2^2}$. For fixed $\lambda >0$, set $v
= \frac{V_+(x)} \lambda$. We claim the following version of Young's
inequality:

\BEL{Young}Uv \leq M(U) + N(v) \EE with equality if and only if $v =
e^U$. To prove this, first consider the case that $v \ge e^U \ge 1$.
Then $N(v) \ge \bar{N} = \int_1^{v} {\rm min}\ \{U(s),\, \ln(s)\}
ds$. The rectangle $[0,U]\times[0,v]$ can be partitioned into two
disjoint regions, with areas $M(U)$ and $\bar{N}$, so $Uv = M(U)
+\bar{N}$. This proves the claim when $v \ge e^U$; the rest is
similar.

By  Green's identity (Lemma \ref{GreenG}), the definition of $U$ and \rq{Young},
\begin{eqnarray}\|\nabla u \|^2_2
&=&
 \int_D |u(x)|^2 V(x) dx
\leq   \int_D |u(x)|^2 V_+(x) dx\nonumber\\
&=& \frac{\|\nabla u \|^2_2}{4\pi}\int_D U(x) V_+(x)dx\nonumber \\
& \leq &   \frac{\|\nabla u \|^2_2}{4\pi} \left( \lambda \int_D
M(U(x))dx + F(\lambda) \right). \label{minF}
\end{eqnarray}

By \rq{MT},  $\dsize \int_D M(U)dx  \le C_2|D|$ and \rq{inf}
follows.

We now show that equality is attained. Flucher proved in \cite{Fl}
that equality occurs in \rq{MT} for some $u_* \in W^{1,2}_0(D)$. We
can assume that $\| \nabla u_*\|_2=1$. Let $U = U_* =4\pi {|u_*(x)|^2}$,
so $\int_D M(U_*)dx = C_2|D|$. By Lemma \ref{EuLaO}, $ -\Delta u_*= Vu_*$, where $V =  \omega^{-1}e^{4\pi|u_*(x)|^2}$ and
$\omega = {\int_D |u_*(x)|^2e^{4\pi| u_*(x) |^2} dx}$. Set $\lambda = \omega^{-1}$, so
$e^{U_*}= V_+ \lambda^{-1} =v$ with equality in \rq{Young} for all $x$. Direct calculation gives $\int_D|U_*V_+|dx
=4\pi$. So, by integrating \rq{Young}
$$4\pi = \lambda \int_D M(U_*(x))\ dx + F(\lambda) = \lambda C_2|D| +  F(\lambda).$$ Thus for
these choices of $u$, $V$ and $\lambda $,  \rq{inf} is an equality.
$\Box$.
\medskip

\section* {3. Minimal support results for other  elliptic equations}
\setcounter{section}{3} \setcounter{Thm}{0} \setcounter{equation}{0}\setcounter{subsection}{0}

\medskip

In this section we prove minimal support results for  other well-known differential equations. Our linear examples are operators in divergence form and Schr\"odinger equations with first order terms. We also study   some related non-linear elliptic equations.
\medskip

\subsection{Operators in divergence form}

\medskip
Our next theorem generalizes  Theorem \ref{ThmX} to operators in
divergence form. Let $a,\ b$ be positive $L^{\infty} (D)$ functions
with $\frac 1a,\ \frac 1b$  in $L^{\infty} (D)$. Define the weighted space $L^{p,b}(D)$ using the norm
$$\| u\|_{p,b}^p=\int_D| u(x)|^p b(x)dx$$
and define $W^{1,2,a}_{0}(D) $ as the closure of $C^\infty_0(D)$
with respect to $\| \nabla u\|_{2,a}$. These norms are equivalent to the ones
with $a\equiv b\equiv 1$, and hence we have the usual compact
embeddings $W^{1,2,a}_0(D)\to L^{2q,b}(D)$ for $q< \ov q$, and for
$q=\infty$ when  $n=1$. When $n>2$ and $q= \ov q$, this embedding is
bounded, but not compact. In what follows, we will denote by
$K=K(D,n,2q,a,b)$ the best constant in the weighted Sobolev
embedding theorem
\begin{equation}
 \| u\|_{2q,b} \leq K  \|\nabla u\|_{2,a}.\label{sobab}
  \end{equation}
Let  $u\in W_0^{1,2,a}(D)$ be a  non-trivial solution  for
  \BEL{divform} -{\rm div} (a \nabla
u )(x)=V(x)b(x) u(x),  \EE  in the sense that
$$\int_D a \nabla u\cdot \nabla \psi\,dx =\int_D  Vu\psi\, b\, dx\quad \forall
\psi \in C_0^{\infty}(D).
$$
\medskip

\noindent
\begin{Thm}\label{Xab1}   Suppose $V\in L^{r, b}(D)$ with $r,\ q$
as in Theorem \ref{ThmX}.  Let     $u\in W_0^{1,2,a}(D)$ be a
non-trivial solution of \rq{divform}. Then
 \BEL{mainab}
K^2 ||V||_{r,b} \ge 1. \EE Equality can occur
when $r>\frac n2$ for $n\geq 2$ and   when $n=1$ and $1<r\le
\infty$.
\end{Thm}

\medskip
\noindent {\it Proof.} Let $V\in L^{r,b}(D)$
with $r\ge \frac n2$.  With Green's identity, (Lemma \ref{GreenG}),
we have the minimal support sequence
\begin{eqnarray}
\| u\|_{2q,b}^2 & \leq & K^2 \int_{D}a| \nabla u|^{2}dx
= K^2 \int_{D}V|u|^{2}b\ dx\label{sob12}\\
& \leq  & K^2 \| V\|_{r,b}\| u^2\|_{q,b}  =
    K^2
\| V\|_{r,b}\| u\|^2_{2q,b},\label{hol*}
\end{eqnarray}
 hence \rq{mainab}.
Now let $r>\frac n2$. Since the weighted norms here are equivalent
to the ones used in Lemmas \ref{Markab} and \ref{EuLa}, the proofs
there still hold; there is a non-negative $u_* \in W^{1,2,a}_0
(D)$ for which $\|u_*\| _{2q,b}=K \|\nabla u_*\| _{2,a}$. It
solves \rq{divform} with $V = c|u_*(x)|^{2q-2}$, where $c= \frac{\|
\nabla u_*\|^2_{2,a}}{\|  u_*\|^{2q}_{2q,b}}.$ Note that \rq{sob12}
is an equation when $u=u_*$. Since $\|  u_*\|^{2q}_{2q,b} = \|  u_*^2 \|^{q}_{q,b} = \|  u_*^{2q-2} \|_{r,b}^r$, we get $$\int_{D}V|u|^{2}b\ dx = c\|  u_*\|^{2q}_{2q,b} = c\| u_*^2 \|_{q,b} \ \|  u_*^{2q-2} \|_{r,b} =
\| u_*^2 \|_{q,b} \ \| V \|_{r,b}
$$ so equality is also attained in   H\"older's  inequality     \rq{hol*}, hence in \rq{mainab} as well. $\Box$

\medskip

\subsection{A result for annuli}

\bigskip

Theorem \ref{ThmX} can be extended  to $V\in L^r(D)$ for $r>1$
(instead of $r>\frac n2$)  in the special case of radial solutions
of the Schr\"odinger equation \rq{Sch} on an annulus  in $\R^n$,
with $n \geq 2$. Let $\rho =|x|$.  Fix $0<c<d< \infty$; let $I = (c,d)$ with the
measure $\rho^{n-1}d\rho$ and let $A=\{ x\in \R^{n} \; :\; c < \rho
< d \}$.   Denote the set of radial functions in $W^{1,2}_{0}(A)$ by
$W^{1,2}_{0,\rm{rad}}(A)$ and define $L^p_{\rm{rad}}(A)$ similarly. If $-\Delta u=Vu$
  holds for $u\in W^{1,2}_{0,\rm{rad}}(A)$, it follows that   $V = -\Delta u /u $   is also a radial function.
  The Laplacian of a radial function
   $u$ is $\Delta u=  \frac{1}{\rho^{n-1}}\frac{\partial}{\partial\rho} \left(\rho^{n-1}
\frac{\partial}{\partial\rho}\right) u$, so $-\Delta u=Vu$ reduces to
   \begin{equation}
    -\frac{\partial}{\partial\rho}\left(\rho^{n-1}
\frac{\partial}{\partial\rho}\right) u=
\rho^{n-1}V(\rho)u\label{temp} \end{equation} with $ u\in
W^{1,2,\rho^{n-1}}_{0} (I)$. Let $\omega_n=|S^{n-1}|$.  Then the mapping
$u(|x|)\rightarrow \omega_n^{\frac 1p} u (\rho)$ gives an isometric
isomorphism $L^p_{\rm{rad}}(A)\rightarrow L^{p, \rho^{n-1}}(I)$. Let
\begin{equation}\label{sobrad} \dsize K_{q, \rm{rad}}(A) =\sup_{ u \in
W^{1,2}_{0, \rm{ rad}}(A)} \frac{||u||_{2q}}{ ||\nabla u ||_2} =
\omega_n^{\frac 1{2q}-\frac 12}\!\!\!\!\!\!\sup _{ u \in W^{1,2, \rho^{n-1}}_{0}(I) }
\frac{||u||_{2q,\rho^{n-1}}}{||  u' ||_{2,\rho^{n-1}}} =
\omega_n^{\frac 1{2q}-\frac 12} K \end{equation}
where $K= K(I, 1, 2q, \rho^{n-1}, \rho^{n-1})$ is  as in \rq{sobab}.
We note that $W^{1,2}_{0,\rm{rad}}(A) \subset W^{1,2}_{0}(A)$
implies $K_{q, \rm{rad}}(A) \leq K_q (A)$.

\noindent
\begin{Thm}\label{annulus}
Suppose the Schr\"odinger equation \rq{Sch} has a nontrivial
solution  $u\in W^{1,2}_{0,\rm{rad}}(A)$   for some  $V\in L^r(A)$, with $ 1 <r \leq \infty $ and $q=r^*$. Then
\BEL{Krad}  K_{q, _{\rm{rad}}}^{2}(A)  ||V_+||_r \ge 1 \EE   and equality  can be  attained.
\end{Thm}

\medskip
\noindent {\it Proof.}  We apply Theorem  \ref{Xab1} with $D= I=(c,d)$, and $a(x)=b(x)=\rho^{n-1}$.
We obtain $ K^2 ||V_+||_{r, \rho^{n-1}}\ge 1.$ By \rq{sobrad}, $K^2  =\omega_n^{1 -\frac 1{q}}K_{q,\rm{rad}}^2(A)$,
and since $$ \omega_n ^{1-\frac 1q} ||V_+||_{r, \rho^{n-1}} =
\omega_n ^{\frac 1r} ||V_+||_{r, \rho^{n-1}} = ||V_+||_{L^r(A)},$$
\rq{Krad} follows.

Let us  show that equality can be attained; by  Lemma \ref{Markab},
there exists     $u_*\in W^{1,2, \rho^{n-1}}_0(I)$   that satisfies
$ \|u_*\| _{2q,\rho^{n-1}}=K\|u_*'\| _{2,\rho^{n-1}}$.
Define $v_*\in W^{1,2}_{0,\rm{rad}}(A)$ by $v_*(x)= u_*(|x|)$.
 We see at once that $\|v_*\| _{2q}=K \omega_n^{\frac{1}{2q}-\
\frac 12} \|\nabla v_*\|_2 =  K_{q,\rm{rad}}(A) \|\nabla v_*\|
_{2}$. So, $v_*$  is an extremal of \rq{sobrad}.
We apply Lemma \ref{EuLa} on $I$, with $a=b=\rho^{n-1}$, to get the
following  equation holding in the distribution sense on $I$:
$-\frac{\partial}{\partial \rho}( \rho^{n-1}
\frac{\partial}{\partial \rho} u_*)= \rho^{n-1} \tilde{V} u_* $,
with
$$\tilde{V}(\rho )= \frac{ |u_*(\rho )|^{2q-2} ||u_*'||_{2,\rho^{n-1}}^2}{\|u_*\|_{2q,\rho^{n-1}}^{2q}} = \frac{ |u_*(\rho )|^{2q-2} ||\nabla v_* ||_{2 }^2}{\|v_*\|_{2q }^{2q}} .$$
Since $v_*(x)=u_*(\rho )$,  by the discussion leading up to
\rq{temp} we get $-{\Delta} v_*=  {V} v_*$, with ${V}(x)=
|v_*(x)|^{2q-2} ||\nabla v_*||_{2}^2\ /\ \|v_*\|_{2q}^{2q}$. Since
$(q-1)r=q$,
$$\left(\int_A |v_*(x)|^{(2q-2)r}dx\right)^{\frac 1r} = \|v_*\|_{2q}^{2q-2}$$
and $$ K_{q, \rm{rad}}^{2}(A) \| V\|_r =  K_{q, \rm{rad}}^{2}(A) \frac {||\nabla v_*||_2^2}{\|v_*\|_{2q}^2 } = 1.\ \
 \Box  $$

\medskip
\medskip

 \subsection{Minimal support results for $-\Delta u= Vu+W \cdot \nabla u$}

\medskip
In this section, we prove  minimal support results for solutions of
second order elliptic equations  with first order terms.
Specifically, we consider

 \begin{equation}
 \label{diffin1} -\Delta u= Vu+W\cdot \nabla u.
 \end{equation}
 Throughout this section, $W$ has values in ${\bf R}^n$  and is defined
 on $D\subset
{\bf R}^n$, with  $n\ge 1$;   $r$ and $q=r^*$ are as in Theorem
\ref{ThmX}. The equation  \rq{diffin1} is assumed to hold in the
distribution  sense, i.e.
$$\int_D\nabla u \cdot \nabla \psi\, dx =\int_D (Vu+W\cdot \nabla
u)\psi\, dx,\quad \forall \psi \in C_0^{\infty}(D).$$

 \begin{Thm}\label{2pot}
Suppose that \rq{diffin1} has a nontrivial
solution $u\in W^{1,2}_0(D)$, with  $V\in L^r(D)$ and  $W\in
W^{1,r}(D, \R^n)$,   $r>  \frac n2$.  Then, \BEL{nec} K_q^2(D)
\left\|V-\frac 12\ {\rm div} W\right\|_r \ge 1 \EE and equality can
be attained.
\end{Thm}

\medskip
\noindent Equality in \rq{nec} can be attained when
 $W\equiv 0$, for then Theorem \ref{2pot} reduces to Theorem
\ref{ThmX}.  The theorem also holds, with the same proof, when $r=  \frac n2$ and $n \ge 3$.  In this case, the inequality  in \rq{nec} is strict, and it is sharp because Theorem \ref{2pot} reduces to Theorem
\ref{ThmX3} when $W\equiv 0$.

\medskip
\noindent {\it Proof.} By taking real or imaginary parts of \rq{diffin1} we can
assume $u$  real-valued.  We assume $n\geq 3$; the proofs for $n=1,\
2$ are similar. By Sobolev's embedding theorem,  $W^{1,r}(D, \R^n)
\subseteq L^{n}(D, \R^n)$, because $r\ge \frac n2 $. Since $|\nabla
u|\in L^2(D)$, H\"older's inequality implies that $\nabla u\cdot W
\in L^{p}(D)$ with $p = \frac{2n}{n+2} = (2\ov q)^*$, also that $V
u\in L^p(D)$.  So  we can apply Green's identity (Lemma
\ref{GreenGG}, with $a\equiv  1$ and $F= Vu+ W\cdot \nabla u$)  to
get

 \BEL{1*} ||u||_{2q}^2  \leq  K_q^2(D) ||\nabla u||_2^2
= K_q^2(D)   \int_D \left(   V u+ \nabla u\cdot W\right)u\,
dx. \EE
  The same argument  used to prove Lemma \ref{GreenG} justifies the identity
\begin{equation}\label{minigreen}
\int_D  2u \nabla u \cdot W dx = \int_D   \nabla (u^2) \cdot W dx = - \int_D u^2\ \div W\ dx.
\end{equation}
From \rq{1*}, \rq{minigreen} and H\"older's inequality  it follows that
\BEL{H} ||u||_{2q}^2 \leq  K_q^2(D) \int_D u^2\left(V -\frac 12\  \div  W\right) dx
\leq K_q^2 (D)||u||_{2q}^2\left\| V-\frac 12\ \div W\right\|_r
.\EE
We conclude that $K_q^2(D)\left\| V-\frac 12\ \div W\right\|_r \ge 1$.  $\Box$

\medskip
In the next theorem, we prove a minimal support result for the
solutions of \rq{diffin1} under weaker assumptions on $W$.

\begin{Thm}\label{2pot1}
Suppose that  the differential equation \rq{diffin1} has a
nontrivial  solution $u\in W^{1,2}_0(D)$ with $V\in L^{r}(D)$ and
$W\in L^s(D, \R^n)$. Let $s\ge 2r \ge n$ (but if $n=2$, let $r> 1$).
Then \BEL{nec2}
 K_q (D)\left( K_q(D)||V||_r + |D|^{\frac 1 {2r} -\frac 1s} ||W||_s \right)\ge 1.
\EE
\end{Thm}

\medskip
\noindent {\it Proof.} Again, we can assume that $u$ is real-valued.
Since $s\ge n$, the proof of \rq{1*} still holds, and gives a
similar formula:
$$||u||_{2q}||\nabla u||_2 \le K_q(D)||\nabla u||_2^2= K_q(D)\left(\int_D   u^2V dx +  \int_D  u\nabla u \cdot W  dx\right)
$$
Applying H\"older's inequality  with exponents $q $
and $ r$ to the first integral, and H\"older's inequality with
exponents $2 q$, $2$, $s$
 and $\left(\frac 1 2-\frac 1{2 q}-\frac 1s\right)^{-1}=\left(\frac 1{2 r} -\frac 1s\right)^{-1} $
 to the second integral (if $s=2r$, the last exponent is not needed),
$$
||u||_{2 q}||\nabla u||_2  \leq K_q(D)
 \left(  ||u||_{2 q}^2 ||V||_r  +  ||u||_{2 q}||\nabla u||_2||W||_s|D|^{\frac 1{2 r} -\frac 1s}\right).
 $$
Applying  Sobolev's inequality \rq{sob}   to the first summand on
the right  hand side,
$$
||u||_{2 q}||\nabla u||_2 \leq  K_q(D)||u||_{2 q}||\nabla u||_2
\left( K_q(D)||V||_r     + |D|^{\frac 1{2 r} -\frac 1s} ||W||_s \right).
$$
So,
$$
1\leq K_q(D)\left(   K_q(D)||V||_r     + |D|^{\frac 1{2 r} -\frac
1s} ||W||_s \right) . \ \ \Box $$

\medskip
We conclude with a corollary to Theorem \ref{Xab1}.
\begin{Thm}\label{exactw}
Suppose that  \rq{diffin1} has a nontrivial solution $u\in
W^{1,2}_0(D)$ for   $V\in L^r(D)$, with $r>\frac n2$, and   for
$W\in L^{\infty}(D, \R^n)$. Suppose also that $W$ is exact,
i.e. $W=\nabla \phi$ for some $\phi \in
W^{1,\infty}( D)$. Let
$$K_{q,\phi}(D)=\sup_{u\neq 0}\frac{\| u\|_{2q,e^{\phi}}}{\|  \nabla u\|_{2,e^{\phi}}}.$$
Then \BEL{neca} K_{q,\phi}^2(D) \left\|V_+\right\|_{r,e^{\phi}} \ge 1 \EE and equality can be attained.
\end{Thm}

\medskip
\noindent
{\it Proof.} Since $-\Delta u-W\cdot \nabla u= Vu$ and $W=\nabla \phi$, we have $-\mbox{div} (e^{\phi}\nabla u )=e^{\phi}Vu.$
Applying  Theorem \ref{Xab1} with $a=b=e^{\phi}$ immediately yields the desired
result.$\Box$
\medskip

Theorem \ref{exactw} has the advantage of being sharp for any exact $W$, but the estimate in Theorem \ref{2pot}
has the advantage that it does not involve the weight $e^{\phi}$.

\bigskip

\subsection{Some nonlinear equations}

\medskip
In this section  we  study  certain  nonlinear differential
equations. We start with the equation \BEL{equbeta} -\Delta u =
V|u|^{\beta-1} u \EE
where $1\leq \beta$.
Here $V$ is assumed to
be real, but $u$ can be complex.  We say that $u\in W^{1,2}_0(D)$
is a {\it very weak solution} of the equation \rq{equbeta}  if
$$\int_D\nabla u \cdot \nabla \psi \,dx= \int_D V|u|^{\beta-1} u\psi\,dx ,\
\forall \psi \in C_0^{\infty}(D).$$
\medskip
\begin{Thm}\label{Sbeta} Assume that \rq{equbeta} has a  nontrivial very weak solution  $u\in W^{1,2}_0(D) $.
Let $\hat q = q(\beta+1)/2$ and assume $\hat q \le \ov q$. If $n\le 2$, let $\hat q<\infty$.
If $V\in L^r(D)$ with $r=q^*$, then
\BEL{betaS} K_{\hat q}^2(D) \|V_+ \|_r\ ||u||_{q(\beta +1)}^{\beta -1} \geq 1.  \EE
Equality can be attained in \rq{betaS} when $\hat q < \ov q$.
\end{Thm}

\medskip
\noindent {\it Proof.} Assume $n\geq 3$; the proof is similar for
$n\le 2$. By   Sobolev's inequality, $u\in L^{2\ov q}(D)$. A
calculation shows that $V|u|^{\beta-1}\in L^{ \ov q ^*}(D)$,
  allowing Green's identity (Lemma \ref{GreenG}A
with $V |u|^{\beta-1}$ replacing $V$) in the minimal support
sequence below.
\begin{eqnarray} \label{nonlinS} ||u||_{q(\beta+1)}^2&=& ||u||_{2\hat q}^2 \leq K_{\hat q}^{2} ||\nabla u||_2^2
 = K_{\hat q}^{2} \int_D |u(x)|^{\beta +1}V(x) \, dx \\ \nonumber &\leq& K_{\hat q}^{2}\int_D |u(x)|^{\beta +1}V_+(x) \, dx
 \le K_{\hat q}^2||u||_{(\beta +1)q}^{\beta +1}||V_+||_{r}
\end{eqnarray}
from which \rq{betaS} follows. If $\hat q < \ov q$, Lemmas \ref{Markab} and \ref{EuLa} provide a $u_* \ge 0$ such that
$ -\Delta u_* =  c u_* ^{2\hat q -1}$ with
$$c = \frac{||\nabla u_*||_2^2}{\|u_*\|_{2\hat q}^{2\hat q}} = K_{\hat q}^{-2}\|u_*\|_{2\hat q}^{2-2\hat q}.  $$
So, $-\Delta u_* = Vu_*^\beta$, which is \rq{equbeta}, with  $V
=V_+= cu_*^{2\hat q  -1 -\beta} = cu_*^{(q-1)(\beta+1)}$. Thus,
$\|V_+\|_r = c\| u_* \|_{q(\beta+1)}^{(q-1)(\beta+1) }=c\| u_*
\|_{2\hat q}^{2\hat q - \beta-1}$, which gives equality in
\rq{betaS}. $\Box$

\bigskip

We now consider  the equation \BEL{egradbeta} -\Delta u = V|\nabla
u|^{\beta}  \EE   with $0< \beta \leq 2$. The case $\beta=2$ is particularly interesting and well studied
in the literature,  (see e.g. \cite{C} and the references cited
there).
We assume that  $u\in W^{1,2}_0(D)$ is a nontrivial weak solution of
\rq{egradbeta},  in the sense that \BEL{Id5} \int_D \nabla u\cdot
\nabla \psi dx = \int_D |\nabla u|^\beta V\, \psi dx \EE for every
$\psi \in W^{1,2}_0(D)$. In the following theorem, we depart from our convention that
$q=r^*$.
\begin{Thm}\label{gradbeta} Let $u$ be as   in \rq{egradbeta} with   $0<\beta<2$,  $V\in L^r(D)$,  $q < \ov q$ and
$\frac{1}{2q} +\frac \beta 2+\frac 1r = 1$.
 Then \BEL{inbeta} K_q(D) ||\nabla
u||_{2}^{\beta-1} ||V||_r \ge 1.
 \EE
When $\beta=2$, $||Vu||_\infty \ge 1$.
\end{Thm}

\medskip
\noindent {\it Proof.}  By Sobolev's inequality and \rq{Id5}, with $
\psi =\ov{u}$, and by H\"older's inequality with exponents $\frac 2
\beta $, $ 2q$, and $r$ we have the following minimal support
sequence
$$||u||_{2q} ||\nabla u||_{2 }  \leq K_q(D) ||\nabla u||_{2 }^2 =
 K_q(D)\int_D    \ov{u} V|\nabla u|^\beta\, dx
 \leq
K_q(D)||u||_{2 q}   ||\nabla u||_{2}^\beta||V||_r,
$$
which implies \rq{inbeta}. When $\beta=2$, \rq {Id5} shows
$\int_D|\nabla u|^2(1-V \ov{u})dx =0$. If $||Vu||_\infty <1$, then
$|1-V \ov{u}|>0$ a.e., so $\nabla u\equiv 0$ on $D$ and $u\equiv 0$, a contradiction.
$\Box$
\medskip
When $n \ne 2$, \rq{inbeta} is also valid for  $q = \ov q$.

\medskip
\section* {4. Minimal support results in one dimension}
\setcounter{section}{4} \setcounter{Thm}{0} \setcounter{equation}{0}\setcounter{subsection}{0}

\bigskip

In this section,  we let $n=1$ and $D=(-b,b)$.  We consider
nontrivial  solutions in the distribution sense
of the equation \BEL{sch1} - u''(x)= V(x) u(x)
\qquad u\in W^{1,2}_0(D). \EE
We can assume without loss of generality that $u$ is real-valued. We can extend $u$ continuously to $[-b,b]$ by setting $u(-b)=u(b)=0$.
As noted elsewhere, most of the results in this paper hold
in this setting, but in this section we show how the case $r=1$
differs.

Thus we consider $V\in L^1(D)$ in Theorem \ref{ThmX}, so
that $q=\infty$. The minimal support sequence still holds in this
case, but the variational work in Lemma \ref{EuLa} does not, so
interesting new questions on sharpness and extremals   arise.
 We prove an analogue  of Theorem \ref{ThmX}, replacing $L^1(D)$ with the space $M$
of signed measures $V$ on $D$, (see e.g. \cite{Ru} for  the
definition and properties of signed measures)  with norm
$$\| V\|_M=|V|(D) < \infty $$
In the  special case where $V\in L^1(D)$, we have $\| V\|_M=|V|(D)=\int_D|V(x)|dx=\| V\|_1.$
\bigskip
\begin{Thm}\label{Steve1}    Assume $u\in W_0^{1,2}(-b,b)$ is a nontrivial solution of \rq{sch1},
with $V\in M$. Then \BEL{rBe1}b||V||_M \ge 2. \EE Equality is
attained when $u=1-\frac{|x|}{b}$. Equality is not possible with $V\in L^1(-b,b)$, but \rq{rBe1}
is still sharp in this case.
 \end{Thm}
\bigskip

\noindent {\it Remark:} For $w\in W^{1,1}(-b,b)$ we have (see \cite{E}, p.286):
\begin{equation}
w(t)-w(s)=\int_s^tw'(\tau)d\tau.\label{ftc}
\end{equation}
We can apply this with $w=u$ and also with $w=u'$; since  $V\in
L^1(-b,b)$,   and $u\in L^\infty  (-b,b)$, \rq{sch1} implies $u''
\in L^1(-b,b)$. We will use \rq{ftc} without further comment
throughout this section.

\medskip
\noindent {\it Proof.} For $x\in (-b,b)$,  by \rq{ftc} and H\"older's
inequality
$$|u(x)| \le \int_{-b}^x |u'| dt \leq \left((x+b)\int_{-b}^x |u'|^2
dt\right)^{\frac 12}.$$ Likewise, $$|u(x)| \leq \left((b-x)\int_x^b
|u'|^2 dt\right)^{\frac 12}. $$ By squaring and algebra,
$$\frac{2|u(x)|^2}b \le \left(\frac{1}{ x+b}+\frac 1{b-x}\right) |u(x)|^2 \le ||u'||_2^2$$   So, $K_\infty^2(-b,b) \le \frac b2$. Now, apply the minimal support sequence for equation \rq{sch1},
 with $q=\infty$. The "H\"older step" in the sequence can be replaced by
 $$ \int u^2\ dV \le ||u^2||_\infty ||V||_M $$ and we get \rq{rBe1}.
\medskip
Setting $b=1$ for simplicity, the claim about $u=1-|x|$
it is easy to verify directly. For this $u$, $\frac V 2$ is a Dirac mass at $x=0$. We see that \rq{rBe1} is sharp for $V\in L^1(-1,1)$ by considering an approximating sequence to $u=1-|x|$. This also implies that
$K_\infty (-b,b)= \sqrt{\frac b2}$.
\medskip
Now, we prove that for $V\in L^1(-b,b)$ equality is never attained,
(this reasoning also gives an independent proof of \rq{rBe1} for
this case). We may assume $u>0$ on $(-b,b)$, for if it changes sign,
we may apply \rq{rBe1} to a restriction of $u$, and we are done. For
now, suppose that $u$ attains it maximum value at a unique point
$c\in (-b,b)$. As in the remark above, $u'$ is defined and
continuous  on $(-b,b)$, and so $u'(c)=0$. Next, we claim that
\BEL{claim} \int_c^b |V|\ dx
>\frac{1}{b-c}.\EE
To prove this, we may assume that $u(c)=1$. By the
mean value theorem on $[c,b]$, there is a point $c<d<b$ such that
$u'(d) = -\frac{1}{b-c}$. Since $u$ is maximal only at $c$, we have
$0<u<1$ and $\frac 1u >1$ on $(c,b)$, and
$$\int_c^d |V|\ dx = \int_c^d \frac{|u''|}u \ dx > \int_c^d |u''| \ dx = \frac{1}{b-c}.$$
The claim follows. Similarly, $(b+c)\int_{-b}^c |V|\ dx
>1$. So, $$||V||_M = ||V||_1 >\frac{1}{b-c}+\frac{1}{b+c} \ge \frac2 b,$$
proving that \rq{rBe1} is strict. We have assumed that $u$ attains a
maximum only at one point $c$; the general case follows by similar
reasoning applied to appropriate restrictions of $u$.
$\Box$

\bigskip

\section* {5. Appendix}
\setcounter{section}{5} \setcounter{Thm}{0} \setcounter{equation}{0}\setcounter{subsection}{0}

In this section, we prove various lemmas needed throughout the paper.
Some already appear in the literature in slightly different form, but are presented
here for completeness. We establish existence of extremals, some variational formulas, and
several versions of Green's identity.

\bigskip
\subsection{Existence of Sobolev extremals }

 We first prove the existence of extremals for the weighted Sobolev inequality  \rq{sobab} used in Theorem \ref{Xab1}.  This applies in other settings, such as Theorem \ref{ThmX}, when the weights are $a=b=1$.
 See \cite{E} for the functional analysis used in the lemmas below.

\begin{Lemma}\label{Markab} Let  $q<\ov q$ for $n\geq 2$,  and $1\leq q\leq \infty$ for $n=1$.
Let $a$, $b\in L^\infty(D)$, with  $\frac 1a$, $\frac 1b\in
L^\infty(D)$;  define $K$ as in \rq{sobab}. There is a nontrivial
and non-negative $u_* \in W^{1,2,a}_0 (D)$  for which \BEL{Inab}
\|\nabla u_*\| _{2,a} =K \|u_*\| _{2q,b}.\EE
\end{Lemma}

\medskip
\noindent
{\it Proof. } Let $B_W$ denote the set of all elements of $W^{1,2,a}_0(D)$ with  $\|\nabla u
\|_{2,a} \leq 1$.  By our assumptions on the weights $a$ and $b$,
the  norms in $L^{q,b}(D)$ and $W^{1,2,a}_0(D)$  are equivalent to the norms with $a\equiv b\equiv 1$.
Thus, $B_W$ is weakly compact in  $W^{1,2,a}_0 (D)$. By the Kondrachov-Rellich Theorem   for $n\geq 2$, and
by  the Arzela-Ascoli theorem for $n=1$,
the inclusion $W^{1,2,a}_{0}(D) \ra L^{2q,b}(D)$ is compact.
Let $\{u_n\}$ be a sequence in $W_0^{1,2,a}(D)$ such that
 $$\lim_{n \rightarrow\infty}\frac{\|u_n\|_{2q,b} }{\|  \nabla
 u_n\|_{2,a}} =K.$$
We can assume by scaling that $\| \nabla u_n\|_{2,a} =1$.
Since  $B_W$ is weakly compact in  $W_0^{1,2,a}(D)$, there exists a
subsequence   $\{ u_{n_k}\} \subset \{u_n\}$  that  converges weakly to
some $u_*\in B_W$. By the compactness of the inclusion $W^{1,2,a}_0(D)\ra L^{2q,b} (D)$, there is a subsequence of    $\{u_{n_k}\}$, that we label again with
 $\{ u_{n_k}\} $,  that converges  to some   $w\in L^{2q,b}(D)$ in the strong topology of $L^{2q,b}(D)$.  That is,
 $\dsize\lim_{k\to\infty} \| u_{n_k}- w\|_{2q,b}=0.$
  But $u_{n_k}\to  w$ also in the weak topology of $ L^{2q,b}(D)$, and so $u^*=w$ a.e.;
 consequently, $w\in B_W$ and $\| \nabla w\|_{2,a} \leq 1$. We have
    \begin{equation}
    K =\lim_{k\to\infty} {\|u_{n_k}\|_{2q,b}}=\|w\|_{2q,b},\label{Ks}
    \end{equation}
 and so
$\dsize \frac{\|w\|_{2q,b} }{\|  \nabla
 w\|_{2,a}}
 \geq K$. But recall that $w=u^*\in W^{1,2, a}_0(D)$, and so
$\dsize \frac{\|w\|_{2q,b} }{\|  \nabla
 w\|_{2,a}}
 \leq K,$
thus proving \rq{Inab}. We can replace $u_*$ by $|u_*|$, if necessary, to get a non-negative extremal, with no effect on \rq{Inab} (see \cite{LL}). $\Box$

\medskip
\noindent

\bigskip
\subsection{Variational work: divergence and Orlicz forms}

We now show that the extremals of the Sobolev inequality
\rq{sobab} solve an equation of the form \rq{divform}, and we give an explicit expression for  $V$
in terms of  these extremals.

\begin{Lemma}\label{EuLa}  Let  $q<\ov q$ with $a$, $b$ and $u_* \ge 0$
  as in Lemma \ref{Markab}.
Then, $-{\rm div}( a \nabla u_*)= b  V u_* $, in the distribution
sense, with \BEL{Vab}V(x)=  \frac{ u_*(x)^{2q-2}  ||\nabla
u_*||_{2,a}^2}{\|u_*\|_{2q,b}^{2q}} .\EE
\end{Lemma}

\medskip
\noindent  {\it Proof.} Let $v(x) = u _*(x) +\epsilon \phi(x)$ where $\phi \in
C_0^{\infty}(D)$ is real-valued.  Let $\delta
=\frac{d}{d\epsilon}|_{\epsilon =0}$. Since $u_*$ is an extremal for \rq{sobab},
 $\dsize \delta \left(\frac{\| u_*+\epsilon \phi \| _{2q,b}^2}{\| \nabla
(u_*+\epsilon \phi)\|_{2,a}^2}\right)=0$.  Direct computation yields
\begin{eqnarray*}
\delta \| u_*+\epsilon\phi\|_{2q,b}^2&=& \frac{1}{q}\|  u_*\|^{2-2q}_{2q,b}\int_D\delta ( b| u_*+\epsilon \phi |^{2q})dx\\
&=&2\| u_*\|_{2q,b}^{2-2q}\int_D u_*^{2q-1}\phi\, b\, dx
\end{eqnarray*}
and
\begin{eqnarray*}
\delta \| \nabla (u_*+\epsilon\phi)\|_{2,a}^2&=&  \int_D\delta (a|\nabla (u_*+\epsilon \phi)|^2)dx\\
&=& 2\int_D (\nabla u_*\cdot \nabla \phi) \ a\ dx.
\end{eqnarray*}
One can justify   passing the derivatives into the integrals by arguing as in \cite{St}. By
the quotient rule
\begin{eqnarray*}
\delta \left(\frac{\| u_*+\epsilon \phi \| _{2q,b}^2}{\| \nabla
(u_*+\epsilon \phi)\|_{2,a}^2}\right)
 &=&
2\frac{||u_*||_{2q,b}^{2 } }{||\nabla u_*||_{2,a}^4}
\int_D \left(\frac{||\nabla u_*||_{2,a}^2}{||u_*||_{2q,b}^{2q }}
u_*^{2q-1}  b\phi -  a(\nabla  u_*\cdot \nabla  \phi )\right) dx =0.
 \end{eqnarray*}
So,  for every $\phi\in C^\infty_0(D)$,  and for  $V$   as
in \rq{Vab}
$$
\int_D V u_* b\,\phi\, dx= \int_D a(\nabla  u_*\cdot \nabla  \phi )  dx
$$
and $-\div( a \nabla u_*)=b Vu_*$, as required. $\Box$
\bigskip

The next lemma shows that the  extremals for the
Moser-Trudinger inequality  in $\R^2$  also satisfy a Schr\"odinger equation.

\begin{Lemma}\label{EuLaO}
Let  $u_*\in W^{1,2}_0(D)$  be a  extremal of the Moser-Trudinger
inequality  \rq{MT}, with  $||\nabla u_*||_2=1$; then $ -\Delta u_*  = V  u_* $ in the distribution sense, where
\BEL{MTEL}V = \frac{ e^{4\pi |  u_*(x) |^2}}{\int_D |u_* |^2\ e^{4\pi|  u_*  |^2} dx.}\EE
\end{Lemma}

\medskip
\noindent {\it Proof.} As in Lemma~\ref{Markab},  we can assume without loss of generality that $u_*$ is non-negative. Let
       $$\dsize U_\epsilon(x)= \frac{4\pi |u_*(x)+\epsilon\phi (x)|^2}{\|\nabla (u_*+\epsilon\phi)\|_2^2} $$
where $\phi \in C_0^{\infty}(D)$ is real-valued. Also set $U=U_0$.  Let $\delta = \frac{d}{d\epsilon} \vert_{\epsilon=0}$. Let $M$ be as in the proof of Theorem  \ref{Orl2}.    Then
     \BEL{MT*} \delta  \left(M(U_{\epsilon})\right)= \int_D e^{U_\epsilon} \delta (U_\epsilon ) dx=0 \EE
     and
$$ \delta (U_\epsilon)   = 8\pi\left (u_*\phi - |u_*|^2 \int_D\nabla u_* \cdot \nabla \phi\, dx\right).$$
One can justify the passing the derivatives into the integrals by arguing as in \cite{St}.
With $\omega=\int_D |u_* |^2\ e^{4\pi|  u_*  |^2}$, we can simplify the equation \rq{MT*}  as follows:
        \begin{eqnarray*}
        0 & =&    \int_D e^{U}  \left(u_*
           \phi    -|u_*|^2 \left(\int_D\nabla u_*\cdot \nabla \phi\, dx\right ) \right ) dx\\
           & =&    \int_D e^{U}   u_*
           \phi  dx  - \omega\int_D\nabla u_*\cdot \nabla \phi\, dx   =  \omega\int_D \left(\omega^{-1} e^{U}  u_*
           \phi  -\nabla u_*\cdot \nabla \phi\right)   dx. \label{eulan}
        \end{eqnarray*}
which proves \rq{MTEL}. $\Box$

\bigskip
\bigskip
\subsection{Green's identities for divergence and Orlicz forms.}

\bigskip
The following lemmas substitute for Green's identity throughout the paper, often  with $a\equiv b \equiv 1$.
\medskip
\begin{Lemma}\label{GreenG}   Let $a(x)>0$ with $a$, $\frac 1a   \in
L^\infty(D)$ and let $b\in L^{\infty}(D)$.  Let $u\in
W^{1,2,a}_0(D)$, with  $-{\rm div}(a\nabla  u)= bV u$ in the distribution sense.
Then the identity
\begin{equation}\label{greenh2}
\int_D  |\nabla  u|^2 adx= \int_D V |u|^2 b dx.
\end{equation}
holds whenever either A),  B) or C) hold:
\medskip

A) $V \in L^{\ov q^*} (D)$ when  $n\ne 2$, and    $V \in L^{r} (D)$
for some $r>1$ when $n=2$,

B) $n\ge 3$ and either $|V(x)| \leq c|x|^{-2}$ or $|V(x)| \leq c
({\rm dist}(x ))^{-2}$ (see \rq{dist}),

C)  $u$ is real-valued,   $a\equiv b\equiv 1$, $n=2$ and   $V\in
L{\rm log}L(D)$.

\end{Lemma}

\medskip
\noindent {\it Proof.} By definition of solution in the distribution sense
\BEL{pde1}   \int_D   \nabla u \cdot \nabla \psi  \, a  dx= \int_D V u \psi b\ dx \EE
for every $\psi \in C^\infty_0(D)$.  The norm in
$W^{1,2,a}_0(D)$  is equivalent to the norm in $W^{1,2}_0(D)$, so
$C^\infty_0(D)$ is dense in both. Let $\{\psi_n\}$ be a sequence of functions in $C^\infty_0(D)$
that converges to $\ov u$ in $W^{1,2,a}_0(D)$. Then
$$
\int_D\nabla  u\cdot  \nabla \psi_n \, a\,  dx -\int_D|\nabla  u|^2 a\, dx \leq ||  \nabla \psi_n - \nabla \ov u||_{2,a} ||\nabla u||_{2,a} \to 0.
$$
To complete  the proof of \rq{greenh2}, it suffices to show that
$Vu\psi_n $ converges  to $V|u|^2$ in $L^{1,b}(D)$ when $V$ is as in
A), B) or C). Assume A), and that $n\neq 2$ (the proof is similar
when $n=2$). By Sobolev's inequality,  $\psi_n$ converges to
$\ov{u}$ in $L^{2\ov{q}} (D)$. By Holder's inequality,
  $$\dsize \int_D |V u \psi_n  -  V |u|^2|\,b\, dx
  \leq  ||b||_\infty ||V  ||_{\ov q*}||u ||_{2\ov{q}}
 \| \psi_n- \ov u\|_{2\ov q}\to 0. $$

For case B), first assume $|V(x)| \leq c|x|^{-2}$. Note that $Vu\psi
\in L^{1,b}(D)$ because $|V|\leq C|x|^{-2}$, and by the
Hardy-Sobolev inequality \rq{H1}, $\psi |x|^{-1}$ and $u |x|^{-1}$
are in $L^2(D)$. Let  $\{\psi_n\}\in C^\infty_0(D)$ be a sequence of
functions that converges to $\ov u$ in $W^{1,2,a}_0(D)$.  We  show
that $V  u \psi_n $ converges in $L^{1,b}(D)$ to $  V |u|^2 $. By
H\"older's inequality and \rq{H1}
$$
\int_D  \left| V  u \psi_n -  V  |u|^2  \right|\, b dx \leq \int_D c
|x|^{-2}| u|\, |\psi_n -\ov u| bdx
 $$
 $$
 \leq  c ||b||_\infty \left(\int_D  |x|^{-2} |u|^2 dx\right)^{\frac 12}
  \left(\int_D |x|^{-2} | \psi_n -\ov u |^2dx\right)^{\frac 12} \leq  c ||b||_\infty
  ||\nabla (\psi_n-\ov u)||_2 \, ||\nabla u||_2 \to 0.
 $$
  This  proves   \rq{greenh2}  in this case. The proof is similar when  $V(x) \leq c({\rm dist}(x))^{-2} $, using the inequality \rq{BM}  instead of \rq{H1}.

Next, we  assume C), and without loss of generality, that $\|\nabla u\|_2= 1$.
Let  $\{\psi_n\}\subset C_0^{\infty}(D)$ converge to $u$ in  $W^{1,2}_0(D)$. We can choose
$\lambda_n \downarrow 0 $ so that  $\| \nabla (u-\psi_n)
\|_2 \lambda_n^{-1} \rightarrow 0$. The Moser-Trudinger inequality
implies
$$\int_D e^{ \left(\frac{u-\psi_n}{\lambda_n}\right)^2} dx\leq \int_D
e^{4\pi\left(\frac{u-\psi_n}{\lambda_n}\right)^2} dx\leq \int_D
e^{4\pi\left(\frac{u-\psi_n}{\|\nabla (u- \psi_n)\|_2}\right)^2}dx<
(C_2+1)|D|<\infty,$$  with  $C_2$ independent of $u$ and $\psi_n$. A similar inequality holds
when $\frac{u-\psi_n}{\lambda_n}$ is replaced by $u$.  Define the
functions $M_2(t) = e^{t} - t - 1$ and $N_2(t) = (t+1)\log(t+1) -
t$  for $t\ge 0$. These are complementary Orlicz functions,
with properties similar to $M$ and $N$, such as Young's inequality \rq{Young} (see also
\cite{KR}).
Using this, the inequality $2|ab|\leq
a^2+b^2$, and H\"older's inequality:
\begin{eqnarray*}
\int_D \left|u \,V (u-\psi_n )\right| dx & = & \lambda_n \int_D
|\frac{u
(u-\psi_n )}{\lambda_n}V| dx \\
& \leq &  \lambda_n \int_D M_2\left(\frac{|u
(u-\psi_n )|}{\lambda_n}\right) dx+ \lambda_n\int_D N_2(|V|) dx \\
 & \leq & \lambda_n \int_D e^{|\frac{u
(u-\psi_n )}{\lambda_n}|}dx + \lambda_n\int_D N_2(|V|) dx\\
 & \leq & \lambda_n \int_D e^{\frac{u^{2}}{2}} e^{\frac{((u-\psi_n)\lambda_n^{-1})^{2}}{2}} dx+ \lambda_n\int_D N_2(|V|) dx\\
 & \leq & \lambda_n\left\{\left(\int_D e^{u^2}dx\right)^{\frac 12}\left(\int_D e^{(\frac{u-\psi_n}{\lambda_n})^2}
 dx\right)^{\frac 1 2}  +\int_D N_2(|V|)dx\right\}.
\end{eqnarray*}
Thus $ u V \psi_n$ converges in $L^1(D)$ to $Vu^2.
\Box$

\bigskip

The next lemma is used in Section 3. It contains Lemma
\ref{GreenG} part A  as a special case.
\begin{Lemma}\label{GreenGG}   Let $a(x)>0$ with $a$, $\frac 1a   \in
L^\infty(D)$.  Let $u\in W^{1,2,a}_0(D)$ be a  solution in the
distribution sense of $-{\rm div}(a\nabla  u)= F$,
where $F \in L^{(2\ov q)^*} (D)$ for $n\neq 2$, and $F \in L^{r}
(D)$ for some $r>1$ for $n=2$. Then,
\begin{equation}\label{gr3}
\int_D  |\nabla  u|^2 adx= \int_D F\, \ov u dx.
\end{equation}
\end{Lemma}

\medskip
\noindent {\it Proof.}  We have $\int_D a \nabla u \cdot \nabla \psi
dx= \int_D F \, \psi \ dx$ for every $\psi \in C^\infty_0(D)$. The
rest is similar to the proof of Lemma \ref{GreenG}, part A.

\end{document}